\documentclass[a4paper,10pt]{article} 
\usepackage[marginparwidth=0pt]{geometry} 
\usepackage{amsmath}
\usepackage{amssymb}
\usepackage{amscd}
\usepackage{authblk}

\usepackage{verbatim}

\title{\bf On humanization of mathematics: aesthetic mathematics\footnote{This is the second version of the previous paper with the same title. This version is a substantial extension of the earlier paper. It includes an important recognition about mathematics: namely, that classical mathematics itself is a form of humanized mathematics. My thoughts on the theme of this paper are still evolving, and I would like to position this paper as a new starting point for my ideas. The development of ideas concerning the foundations and philosophy of mathematics, including technical results, has been proceeding vividly and rapidly in recent years. However, from my perspective, there are still new aspects to be considered. Therefore, I present my thoughts here as a preliminary report.}}

\bigskip

\author{Takao Inou\'{e}}

\affil{Faculty of Informatics,Yamato University, \\
Katayama-cho 2-5-1, Suita, \\
564-0082, Osaka, Japan\footnote{ 
inoue.takao@yamato-u.ac.jp; (Personal) takaoapple@gmail.com \\
I prefer my personal e-mail.} \\
$\enspace$\\
\bf The second version\rm
}

\date{January 10, 2025}

\begin{document}

\maketitle

\begin{abstract}
This paper examines various methods and ideas for humanizing mathematics. The term "humanizing mathematics"—which includes elements of "aesthetic mathematics"—refers to approaches that emphasize the aesthetic, philosophical, and subjective dimensions of mathematics. These approaches aim to make mathematics humanize. It proposes novel directions in mathematics, stressing that the ideas presented here are provisional. Furthermore, it argues that mathematics can take multiple humanized forms. Even traditional mathematics can be interpreted as a form of humanized mathematics. To support this perspective, several mathematical observations are provided. Additionally, a general approach to addressing the Riemann Hypothesis-focusing on proof by contradiction and mathematical logic, particularly model theory—is outlined. The paper also briefly reflects on my prior research through this idea on humanization and concludes with general remarks and future directions.

\medskip
SMC2020: 00A30, 00A05, 00A30, 03A05

%
%
\end{abstract}

\smallskip


\tableofcontents


\section{Introduction}
This paper examines various methods and ideas for humanizing mathematics. The term "humanizing mathematics"—which includes elements of "aesthetic mathematics"—refers to approaches that emphasize the aesthetic, philosophical, and subjective dimensions of mathematics. When discussing "aesthetic mathematics," this paper highlights the aesthetic aspects within this perspective. These approaches aim to make mathematics humanize. It proposes novel directions in mathematics, stressing that the ideas presented here are provisional. Furthermore, it argues that mathematics can take multiple humanized forms. Even traditional mathematics can be interpreted as a form of humanized mathematics. To support this perspective, several mathematical observations are provided. Additionally, a general approach to addressing the Riemann Hypothesis-focusing on proof by contradiction and mathematical logic, particularly model theory—is outlined. The paper also briefly reflects on my prior research through this idea on humanization and concludes with general remarks. 

The structure of this paper is as follows. The ideas are categorized into three distinct eras—the white era, the blue era, and the green era—based on their conceptual evolution over time.

Sections 2 and 3, which form the core of this paper, provide mathematical remarks related to this perspective. Section 2 introduces the white era, focusing on critiques and remarks on traditional mathematics. Section 3 explores the blue era, where it is argued, among other points, that traditional mathematics can be interpreted as a form of humanized mathematics. In addition, a general approach is proposed to address the Riemann Hypothesis, including methods for its proof or disproof. In addition, in Section 3, I propose modal mathematics and modal monadology as new directions for research.

Section 4 introduces the green era, presenting a forward-looking argument for future directions.

Section 5 reflects briefly on prior research from this perspective. Section 6 broadens the discussion by exploring the humanization of mathematics from a multidisciplinary viewpoint across various academic disciplines. Finally, Section 7 concludes with general remarks and future directions.

%
%
%
%
%
%
%
%
%

\section{Mathematical remarks, Part I: my white era till 2023}
First I would like to tell the readers that I am \it a traditional mathematician \rm as a mathematician. I am no math maverick.

\medskip

It seems to me that almost all the mathematics are satisfied with the current way to define mathematical objects to study: mathematical objects are mostly defined by their properties which are formulated in the form of formulas, axioms or axiom-schemata, assuming implicit use of classical predicate logic with equality in bothe object- and meta-levels, i.e. formal syntactic treatment in mathematics and informal thinking in mind. Does nobody have a doubt about  this situation? I know very well that mathematicians prefer to take more and more useful form of definition in a practical sense to conceptually or aesthetically better one, though they might have thier own aesthetics in fact. Sometimes I however do not find a current form of  definitions so savory. This discontent of mine about definitions may gradually be resolved by considering some other form of definitions and mathematics. Further by considering mathematic as something more human and aesthetic. I shall show some kinds of believe about mathematics. That is the following.
\medskip
 
1. Mathematics is very personal. Thus mathematics may take a variety of forms, depending on each indivisual. This idea makes me a bit recall Gottfried Wilhelm Leibniz, Srinivasa Ramanujan, Luitzen Egbertus Jan Brouwer and Florentin Smarandache.
\medskip
 
2. A mathematical proof is a story like a novel.
\medskip

3. I require mathematics esthetic tastes, e.g. in a figuary sense like literature one, musical one and etc.
\medskip

4. I like L. E. J. Brouwer's idea for mathematics: mathematical objects should be constructive. 
\medskip

5. I am not satisfied with only one mathematical fact. For example, further thinking of many embeddings is possible, when one knows one well-known embedding in mathematics.
\medskip

6. Natural numbers are not created by the God. Natural numbers are created by human. Kronecker is not right.
\medskip

7. 1  is very difficult to understand as a mathematical object. I think that the recognition of 1 still has a lot of secrets for mathematics. I think that this may be connected to the resolution of Riemann hypothesis.
\medskip

8.  0 is not a discovered number. 0 is a number created by human.
\medskip

9. Right understanding of right proof is like a sky without clouds.
\medskip

10. Thinking of modality is that of higher structure. On the other hands, thinking of modality is also that of restricting usual thinking of mathematics. 
\medskip

11. The pair of the set of possible worlds and the set of accessiblity relations may become parameter in mathematics, physics, even science and philosophy (Inou\'{e} \cite{Inoue2023a}).
\medskip

12. There would be humanized mathematics or aesthetic mathematics. One example is epistemic mathematics and constructive one as there exist at present. However epistemic mathematics and constructive one are only parts of humanized mathematics or aesthetic mathematics I intend.
\medskip

13. Mathematics is not only one of Greek and European thoughts but one of human thoughts as a variety of music and scales exists in the world.
\medskip

14. Good proofs are very natural. 
\medskip

15. Past approaches to Riemann hypothesis would be not natural if one accept my idea that good proofs are very natural.  If Riemann hypothesis is solved, its proof should be very natural.
\medskip

What I can now do for finding a way to humanize mathematics is to consider and to develop the above 15 remarks.

\section{Mathematical remarks, Part II: my blue era from 2023 to 2025, A short but condensed era of thinking 
}
In this section, among others I will present some novel directions in mathematics. I believe there are many forms of humanized mathematics. Even traditional mathematics, which we call classical mathematics, is itself a form of humanized mathematics. This part is the continuation and developed thought from my white era.
\medskip

1. The term "humanizing mathematics"—which includes elements of "aesthetic mathematics"—refers to approaches that emphasize the aesthetic, philosophical, and subjective dimensions of mathematics. 

2. Mathematical logic consists of four disciplines: proof theory, model theory, recursion theory and set theory. These are individually very deep theories. Most mathematicians, including myself, recognize that modern mathematics is fundamentally based on set theory, particularly axiomatic frameworks such as Zermelo-Fraenkel (ZF) set theory. They more or less know that set theory is very important as knowledge. However, I dare to point out that they do not really know the importance and significance of mathematical logic. It would partly depend on the education which they received in the past. Usual mathematicians (especially in Japan) tend to think that mathematical logic is a kind of philosophy, not mathematics. This recognition is totally wrong. Mathematical logic is the very mathematics. Mathematical logic contains something humanized and modal concepts to be studied. Mathematical logic has a treasure box of humanization and mathematical materials to be humanized. 

3. Negation has been as part of the foundation of mathematics. Griss even proposed negationless mathematics (e.g. Griss \cite{Griss1946}). \footnote{According to Ferguson, D. Nelson's strong constructible negation is compatible with Griss's principles (Ferguson \cite{Ferguson2023}).} Others have tried to construct the meaning of negation in syntactic form, e.g. Fitch's Basic Logic (cf. Fitch \cite{Fitch1942} and Inou\'{e} \cite{Inoue1984}). The treatment of negation reflects the symmetry of the theory. Classical mathematics is symmetric in logic. On the other hand, for example, intuitionistic mathematics is asymmetric in logic. Each theory has a form of humanized mathematics. Paraconsistent logic is not symmetric in the sense of negation, because it contains some contradictions. So paraconsistent theory has an aspect of humanized mathematics. Negation would be something to be constructed according to its meanings.  

4. Symmetry or commutativity represent a form of humanized mathematics. Also so do asymmetry or noncommutativity.

5. I do not think that the inferences in classical mathematics are right. Those are creations in Greek people, especially Euclid. Traditional mathematics is a kind of human creation (refer to Bourbaki \cite{Bourbaki1994}). So I think that traditional mathematics, usually called classical mathematics, one form of humanized mathematics. Of course I do not negate the significance of classical mathematics, thinking of the usefulness in Physics.

6. Fuzzy theory is a kind of humanized mathematics. There is the active research area for  fuzzy intuitionistic set theory.

7. It seems to me that contemporary theories tend to noncommutative ones. This means humanization of mathematics occurs in many academic area.

8. It is remarkable that the Riemann Hypothesis remains unsolved. If the Riemann Hypothesis were to be resolved, the method employed would likely involve proof by contradiction. Should traditional mathematical approaches prove inadequate, mathematical logic—particularly model theory—might offer a pathway to establishing the validity of the hypothesis. Furthermore, my modal mathematics could contribute to this development.

9. Non-standard mathematics is just a usual mathemtaics. It is not a special one as Halmos had thought so. For the references of non-standard mathematics, see the reference of this paper.

10. Noncommutativity is one key notion of mathematics. The strategy to take inverse in making new theory is very important. Recently there is an academic report that, in biology, taking inverse strucure, i.e. kimera creates new creature.

11. Already existing mathematics contains the elements of humanization, e.g. simplification of axioms (axiom schemata) etc.

12. There are many thoughts about mathematics and the foundation which should be considered, e.g. homotopy type theory, topos theory and so on.

13. Category theory is a kind of humanized mathematics using arrows and objects.

14. Scheme theory in algebraic geometry is a kind of humanized mathematics by Alexander Grothendieck (\cite{Grothendieck}).

15. My concept of modal mathematics represents a form of humanized mathematics (see Inou\'{e} \cite{Inoue2025}). The outline of this framework is as follows: mathematical objects are interpreted in the context of modal operators and possible world semantics. Specifically, for any mathematical object $a$ let its name be $aN$. There exists a modal operator $\Box_{aN}$ and the object $\Box_{aN} a$. $a$ is identified with $\Box_{aN} a$. This object is further interpreted through the model associated with possible world semantics for $\Box_{aN} a$. The coordinate system used aligns with the fifth-dimensional framework described in my earlier work (Inou\'{e} \cite{Inoue2023a}). The theory is developed within the framework of a specific epistemic formal system, representing a highly generalized epistemic mathematics. This approach is inspired by Gottfried Wilhelm Leibniz's Monadology. Furthermore, musical scores could be represented within the framework of modal mathematics.

16. Model theory for $\lambda$-calculus created by Dana Scott is a kind of humanized mathematics (Scott \cite{Scott}).

17. Consider more predicative analysis and impredicative one (refer to Pohlers \cite{Pohlers2009}). 

18. Implication would be a transformation from some world to some world.

19. In order to understand and interpret Leibniz's monadology (\cite{Leibniz}), say $\bf LMo$, I shall propose the modal monadology. Monads can be interpreted as follows. For any monad $a$, say the name of monad is $aN$, then there is a modal operator $\Box_{aN}$ and the monad $a$ is identified with $\Box_{aN} a$ and interpreted by the model associated with possible world semantics for $\Box_{aN}$. This is my theory: the modal monadology, say $\bf MMo$. This idea was derived from my modal mathematics.  I think my idea would be inspired by the monadology by Gottfried Wilhelm Leibniz. Clearly we have 
$$\bf LMo \cong MMo$$
That is, $\bf LMo$ is isomorphic to $\bf MMo$. For the detail, refer to my future paper Inou\'{e} \cite{Inoue2025monad}.

20. Topology would govern the all mathematics.

21. Priest's characterization of logic as applied mathematics (\cite{Priest}) is understandable. However, from my perspective, it lacks a certain appeal. One might simply state that logic is a branch of mathematics.

22. To verify the correctness of mathematical theorems using computers, various proof assistants, such as Mizar, HOL, and Coq, have been developed. However, the foundational bases of these systems differ significantly. For instance, Mizar is grounded in Tarski-Grothendieck set theory combined with classical first-order predicate logic. Similarly, HOL is based on higher-order logic, while Coq is built upon intuitionistic type theory. This diversity suggests that mathematics possesses an aesthetic dimension, reflected in the variety of its representational frameworks, at least from my perspective.

23. Characterization theorems in mathematics serve as representatives of the aesthetic aspects of the discipline.

24. The humanization of mathematics can serve as an inspiration for mathematical research.

25. My explanation of humanized mathematics and aesthetic mathematics is currently insufficient. I intend to elaborate on this topic in greater depth in future work.

\section{Mathematical remarks, Part III: my green era}

I do not know what the future holds. As my thoughts on mathematics and the philosophy of mathematics continue to develop, I envision that some of the descriptions of mathematics I propose may come to fruition in what I refer to as my "green era." In that era, the humanization of mathematics might be regarded as a natural aspect of the essence of mathematics itself. This would align with the idea that mathematics is an innate aspect of human cognition, much like Chomsky's assertion that language is innate to humans (Chomsky \cite{Ch1965}).

\section{An explanation of my papers from the viewpoint of this essay}

I will provide brief explanations of my previous papers from the perspective of this essay. My intention is not to promote my earlier works but to discuss them because they are directly relevant to the theme of this paper.

\begin{quote}
[1] Takao Inou\'{e}\rm , Fitch's basic logic as a point of departure -the philosophical foundations of constructive mathematics- (in Japanese),  Kagaku Tetsugaku, Vol. 17, pp. 133-148, 1984, Waseda University Press, 1984.
\end{quote}

I proposed that the foundation of constructive mathematics consists in Fitch's basic logic.  In this system, the nagation is construced not as primitively given object.  This paper is  related to 3 in \S 3.

\begin{quote}
[2] Takao Inoué, A note on unprovability-preserving sound 
translations, Logique et Analyse (N.S.), Vol. 33, pp. 243-257, 
1990.
\end{quote}

This paper is concerned with the above 5 in \S 2 in the previous section.

\begin{quote}
[3] Takao Inoué, A definition of field from a philosophical point 
of view, Abstracts of Papers Presented to the American 
Mathematical Society, Vol. 13, p. 337, 1992.
\end{quote}

This paper is concerned with the number of axioms of fields. There are some similar works about Boolean algebras, etc. 

\begin{quote}
[4] Takao Inoué, Some type of formulas and the faithfulness of 
Flagg and Friedman's translation, Bulletin of the Section of Logic 
(Lodz, Poland), Vol. 21, pp. 2-11, 1992.
\end{quote}

This paper is to find some set of S4-formulas such that the unprovability of the formula of the set preserves in a special case of Flagg and Friedman's translation. I was not satisfied with the only fact that Flagg and Friedman's translation is not faithful.

\begin{quote}
[5] Takao Inoué, Corrections and additions to my paper "A note on
unprovability-preserving sound translations", more general 
constructions, Logique et Analyse (N.S.), Vol. 39, pp. 335-367, 
1996.
\end{quote}

This paper is also concerned with the above 5 in \S 2.

\begin{quote}
[6] Takao Inoué, How to Make a New Logic, arXiv:2108.05934 
[math.LO], 2021, Cornell University.
\end{quote}

This paper is very suggestive one to make a new logic, takeing the inverse arrow system in a given Gentzen system.

\begin{quote}
[7] Takao Inoué, A sound interpretation of Lesniewski's epsilon in 
modal logic KTB, Bulletin of the Section of Logic, Vol. 50, No. 4, 
pp. 455-463, 2021.
\end{quote}

This paper gives a conjecture that Lesniewski's propositional ontology can be embedded in modal logic KTB. The soundness of it is proved. The proof is very natural. In the sense, it is also related to the 14 in \S 2. So is to the above 1, 5, 10  in \S 2.

\begin{quote}
[8] Takao Inoué, On Blass translation for Lesniewski's 
propositional ontologyand modal logics, Studia Logica, 
Vol. 110, No. 1, pp. 265-289, 2022.
\end{quote}

This paper is related to existence, deontic concepts, provability, normal modality. This is also concerned with the above 10 in \S 2. That is, the relation between modality and second order logic.

\begin{quote}
[9] Takao Inoué, On the fifth dimension as the set of possible 
worlds with treestructures. Center For Computer Science Mathematics And Engineering Physics, publications, 2024.
\begin{verbatim}
https://www.researchgate.net/publication/

368155494_On_the_fifth_dimension_as_the_

set_of_possible_worlds_with_tree_structures
\end{verbatim}
\end{quote}

This paper is concerned with the above 11 in \S 2.

\begin{quote}
[10] Takao Inoué, Epistemic systems and Flagg and Friedman's translation. Center For Computer Science Mathematics And Engineering Physics, publications, 2024.
Also arXiv:2307.02688 [math.LO], 2024, Cornell University.
\end{quote}

This paper treats the inverse flow of thoughts in metamathematics.

\begin{quote}
[11] Takao Inoué and Tadayoshi Miwa, Nontrivial single axiom schemata and their quasi-nontriviality of Leśniewski-Ishimoto's propositional ontology $\bf L_1$. arXiv:2402.07030 [math.LO], 2024, Cornell University
\end{quote}

This is a novel presentation of making axiom schemata to be single one.

\section{The humanization of mathematics in a broader context}
A key consideration is the increasing importance of connecting academic disciplines to humanity. This trend is evident not only in the natural sciences but also in fields such as business informatics and marketing within the social sciences. The process of humanization is already underway and should be further advanced across a broad range of disciplines.

In mathematics, the interplay between human creativity and machine computation demands deeper exploration. This paper emphasizes the need to investigate this relationship to guide future academic and practical developments.

Humanization is an ongoing process that must continue across diverse academic fields. In the context of the current AI era, the relationship between AI and humanity—particularly in areas such as AI ethics—has often been neglected and requires urgent attention. I stress the necessity of addressing this connection, as its current treatment remains inadequate.

\section{Some general remarks and future directions}
I would like to repeat the content of the footnote for the title of this paper. This is the second version of the previous paper with the same title which appeared in arXiv as a preprint. This version is a substantial extension of the earlier paper. It includes an important recognition about mathematics: namely, that classical mathematics itself is a form of humanized mathematics. My thoughts on the theme of this paper are still evolving, and I would like to position this paper as a new starting point for my ideas. The development of ideas concerning the foundations and philosophy of mathematics, including technical results, has been proceeding vividly and rapidly in recent years. However, from my perspective, there are still new aspects to be considered. Therefore, I present my thoughts here as a preliminary report.

From my perspective, I would like to repeat that one significant point is that the connection to humanity is becoming increasingly important not only in the natural sciences but also across many fields, including the social sciences, such as business informatics and AI ethics.

My future work will focus on continuing to humanize traditional mathematics and obtaining new results inspired by this humanization and related reflections.

While I have presented some speculative ideas in \S 4, these remain uncertain.

\medskip

Takao Inou\'{e}

Faculty of Informatics,Yamato University, \\
\indent Katayama-cho 2-5-1, Suita, \\
\indent 564-0082, Osaka, Japan

\bigskip 

inoue.takao@yamato-u.ac.jp

(Personal) takaoapple@gmail.com

I prefer my personal email.

\end{document}